\documentclass[11pt,reqno]{amsart} 

\usepackage{amssymb,latexsym}
\usepackage{cite} 

\usepackage{fontenc}
\usepackage{times}

\usepackage{amsfonts,amssymb,epsfig}
\usepackage{amscd}
\usepackage{color}
 \usepackage{hyperref}

\usepackage{textcomp}
\usepackage{mathrsfs}

\usepackage[height=200mm,width=140mm]{geometry} 

\theoremstyle{plain}

\theoremstyle{remark}

\theoremstyle{example}

\numberwithin{equation}{section} 

\theoremstyle{definition}
\newtheorem{definition}{Definition}[section]

\newtheorem{theorem}[definition]{Theorem}
\newtheorem{proposition}[definition]{Proposition}

\newtheorem{lemma}[definition]{Lemma}

\newenvironment{prof7}[1][Proof]{\textbf{#1.} }{\ \rule{0.5em}{0.5em}}

\begin{document}
\title[Induced Semi Riemannian structures on null submanifolds]{Induced Semi Riemannian structures on null submanifolds}
\author{KARIMUMURYANGO M\'en\'edore}

\address{Applied Statistics Institute\\
	university of Burundi\\ 
	\noindent
	P.O.Box:5158 Bujumbura-Burundi 
}
\email{kmenedore@gmail.com}

\author{NDAYIRUKIYE Domitien}
\address{Ecole Normale Sup\'erieure\\ 
	Section de Math\'ematiques\\
	P.~O. Box 6983 \\ 
	Bujumbura\\ Burundi}
\email{domitienndayi@yahoo.fr}

\author{NIBARUTA Gilbert}
\address{Ecole Normale Sup\'erieure\\ 
	Section de Math\'ematiques\\
	P.~O. Box 6983 \\ 
	Bujumbura\\ Burundi}
\email{nibarutag@gmail.com}
\author{NIBIRANTIZA Aboubacar}
\address{Universit\'e du Burundi\\ 
	Institut de P\'edagogie Appliqu\'ee\\
	D\'epartement de Math\'ematiques\\
	Bujumbura-Burundi}
\email{aboubacar.nibirantiza@edu.bi}

\begin{abstract}
	In this paper, we induce a semi-Riemannian metric on the $r-$null submanifold. 
	We establish the links between the null geometry and basics invariants of the associated 
	semi-Riemannian geometry on  $r$-null submanifold and semi-Riemannian constructed from 
	a semi-Riemannian ambient.
\end{abstract} 


\maketitle

\section{Introduction}\label{Section1}
One of the most fundamental problems in submanifold theory is the problem
of isometric immersibility. The embedding problem had been around since Riemann in $1854$.
Soon after Riemann introduced the notion of a manifold, Schl\"{a}fli conjectured that
every Riemannian  manifold
could be locally considered as a submanifold of an Euclidean space with dimension
$\frac{1}{2}n(n+1)$ in $(1873)$. This was later proved in different steps by Janet $(1926)$,
E. Cartan revised Janet's paper with the same title in $(1927)$.
This result of
Cartan-Janet implies that every Einstein $n-$manifold $(n \geq3)$ can be locally
isometrically embedded in an Euclidean $\frac{n(n+1)}{2}$-space
$\mathbb{E}^\frac{n(n+1)}{2}$. 
In $1956$ J.Nash proved that every closed Riemannian $n-$manifold can be isometrically 
embedded in an Euclidean 
$m-$space $\mathbb{E}^m$ with $m =\frac{1}{2}n(3n + 11)$ and he proved also that 
every non-closed Riemannian n-manifold can be isometrically embedded in
$\mathbb{E}^m$ with $m = \frac{1}{2} n(n + 1)(3n + 11)$.
R.E.Green improved Nash's result in $(1970)$  and proved that
every non-compact Riemannian n-manifold can be isometrically embedded
in the Euclidean $m-$space $\mathbb{E}^m$ with $m = 2(2n+1)(3n+7)$. 
In $1965$, A.Friedman proved that any $n-$dimensional semi-Riemannian manifold of index 
$q$ 
can be isometrically embedded in a semi-Euclidean 
space of dimension $\frac{1}{2}  n (n + 1)$ and index $\geq q$. 

In $1970$ Clarke and Green proved that any semi-Riemannian
$n-$manifold $M^n_q$ with index $q$ can be isometrically embedded in a semi-
Euclidean $m-$space $\mathbb{E}_{\overline{q}}^m$ , for $m$ and ${\overline{q}}$ large enough.
Moreover, this embedding may be taken inside any given open set in 
$\mathbb{E}_{\overline{q}}^m$. The problem of discovering simple sharp relationships between intrinsic and extrinsic invariants of a Riemannian submanifold
becomes one of the most fundamental problems in submanifold theory. The main extrinsic invariant is the squared mean
curvature and the main intrinsic invariants include the classical curvature invariants namely the scalar curvature and the
Ricci curvature.

In \cite{Karimumuryango}, the author immersed a null  hypersurface equipped with the Riemannian metric 
(induced on it by the rigging) into
a Riemannian manifold suitably constructed on the Lorentzian manifold and she established
the basic relationships between the main extrinsic invariants and the main intrinsic invariants 
named Chen-Ricci inequality of the null hypersurface in the Lorentzian manifold.
This rigging technique has been introduced 
in\cite{Gutierrez}.

In the present paper, we induce a semi-Riemannian metric on the null submanifold and we immersed isometrically the null 
submanifold endowed with semi-Riemmannian metric in semi-Riemmannian manifold.
We then establish links between the null geometry and basics invariants of the associated 
nondegenerate geometry such as linear connection, the curvature tensor, Ricci curvature
such that is symmetric and  sectional  curvature. 
The remaining of this paper is organized as follows. 
Section (\ref{Section2}) contains most of the prerequisites material of null submanifolds 
in semi-Riemannian Manifold.
The normalization and the associated semi-Riemannian structure on a normalized null 
submanifold using Rigging techniques are introduced and discussed in Section 
(\ref{Section3}).  
The relashionship between the null and the associated semi-Riemannian geometry is considered 
in Section (\ref{Section4}).
  
 \section{Preliminaries on null submanifolds}\label{Section2}
 
 Let $(\overline{M}, \overline{g})$ be a real $(n+k)$-dimensional semi-Riemannian manifold of constant
 index $q\in\{1,\cdots, n+k-1\}$ where $n>1,k\geq 1$. Suppose $M$ is a $n$-dimensional 
 submanifold of $\overline{M}$.  
 In case $\overline{g}_x$ is non-degenerate  on  $T_x M$, then $T_xM$ and $T_xM^\bot$
 are complementary orthogonal vector subspaces of $T_x\overline{M}$, 
 a part of the normal vector bundle $TM^\bot$ (the radical distribution) lies in the tangent 
 bundle $TM$ of a submanifold $M$ of a semi-Riemannian manifold $\overline{M}$.
 Otherwise,$T_xM$ and $T_xM^\bot$ are degenerate orthogonal subspaces but no longer 
 complementary subspaces, that is $TM\cap TM^\bot\neq \{0\}.$ Thus a basic problem of the null  submanifolds is to replace
 the intersecting part by a  vector subbundle whose sections are nowhere  tangent to $M$. 
 To overcome with this problem posed by null submanifolds, the authors  
 Bejancu and Duggal introduced the notion
 of screen distribution  which provides a direct sum decomposition of $T\overline{M}$
 with certain nice properties \cite{Duggal}. Used a screen distribution $\mathcal{S}(TM)$ on $M$ 
 and a screen vector bundle $\mathcal{S}(TM^\bot)$ over $M$ to construct  a transversal
 bundle $tr(TM)$. They obtained the structure equations of $M$ that relate the curvature 
 tensor of $\overline{M}$ with 
 the curvature  tensor of the linear connections induced on the vector bundles involved
 in the study.
 If the mapping $Rad(TM)$ 
 is a smooth distribution with constant  rank $r>0$,
 then, it is said the radical (null, lightlike) distribution on $M$. Also,
 $g$ is called $r-$null ($r-$lightlike,$r-$degenerate) metric on $M$ the submanifold $M$ 
 is said to be $r-$lightlike ($r-$null, $r-$degenerate)
 submanifold of $\overline{M}$, with nullity degree $r$ and  is simply called
 null (lightlike) submanifold.  Any complementary (and hence orthogonal) distribution $\mathcal{S}(TM)$ of $Rad(TM)$ 
 in $TM$ is called a screen distribution.
 For a fixed screen distribution $\mathcal{S}(TM)$ on $M$, the tangent bundle splits as 
 \begin{equation}\label{Radical}
 T M = RadT M \oplus_{orth} \mathcal{S}(T M).
 \end{equation}
 Certainly, $ \mathcal{S}(T M )$ is not unique, however it is canonically isomorphic 
 to the factor
 vector bundle $T M^{*} = T M/ Rad T M$ considered by Kupeli \cite{Kupeli}.
 $ \mathcal{S}(T M )$ can be constructed by using the local equations of the submanifold and 
 therefore it enables us to obtain the main induced geometrical objects: induced connection, 
 second fundamental form, shap operator, $\cdots$
 A screen transversal vector bundle $\mathcal{S}(TM^\bot)$ on $M$ is any
 (semi-Riemannian)
 complementary vector bundle of 
 $Rad(T M)$ in $TM^{\bot}$. It is obvious that  $\mathcal{S}(TM^\bot)$ 
 is non-degenerate
 with respect to  $\overline{g}$ and $TM^\bot$ has the following orthogonal direct 
 decomposition
 \begin{equation}\label{nullscreen}
 TM^\bot=Rad TM\oplus_{orth} \mathcal{S}(TM^\bot).
 \end{equation}

 $\mathcal{S}(TM)$ and $\mathcal{S}(TM^\bot)$ are called a screen distibution
 and a screen transversal vector bundle of $M$ respectively.    
 As $\mathcal{S}(TM)$ is not degenerate, let $\mathcal{S}(TM)^\bot$ be its
 complementary orthogonal in $T\overline{M}\arrowvert_M.$
 Then we have the following decomposition
 \begin{equation}\label{I3}
 T\overline{M}\arrowvert_M=\mathcal{S}(TM)\oplus_{orth} \mathcal{S}(TM)^\bot.
 \end{equation}
 Note that $\mathcal{S}(TM^\bot)$ is a vector subbundle of $\mathcal{S}(TM)^\bot$ 
 and since both are non-degenerate we have the following  orthogonal direct decomposition
 
 $$\mathcal{S}(TM)^\bot=\mathcal{S}(TM^\bot)\oplus_{orth} \mathcal{S}(TM^\bot)^\bot.$$
 
 Since the theory of null submanifold M is mainly based on both $\mathcal{S}(T M )$ and
 $\mathcal{S}(T M^\bot)$, a null submanifold is denoted by $(M, g, \mathcal{S}(T M ), 
 \mathcal{S}(T M ^\bot))$ and we have four sub-cases with respect to the dimension and 
 codimension of $M$ and rank of $Rad TM$:
 
 \begin{itemize}
 	\item [(a)] $r-$null(lightlike) if $r<min(n, k)$;
 	\item [(b)] coisotropic if $r=k < n$   (hence $\mathcal{S}(TM^\perp)=\{0\}$);
 	\item [(c)] isotropic if $r=n<k,$ (hence $\mathcal{S}(TM)=\{0\});$
 	\item [(d)] totally null if $r=n=k$, (hence $\mathcal{S}(TM)=\{0\}
 	=\mathcal{S}(TM^\perp)$).
 \end{itemize}
 We have seen from the above that the normal bundle $TM^\bot$ is orthogonal 
 but not a complement to $T M$, since 
 it intersects the null tangent bundle $Rad TM$.  This creates a problem as 
 a vector of $T_x\overline{M}$ cannot be decomposed 
 uniquely into a component tangent of $T_x M$ and a component of $T_x M^\bot$ .
 Therefore, the standard 
 definition of second fundamental forms and the Gauss-Weingarten formulae do not work.
 To deal with this anomaly, one of the techniques 
 used consists to split the 
 tangent bundle $T\overline{M}$ into four non-intersecting complementary 
 (but not orthogonal) vectors bundle.
 
 \begin{theorem}\cite{Duggal}
 	Let $\big(M, g,\mathcal{S}(T M ), \mathcal{S}(T M ^\bot)\big)$ be an $r-$lightlike 
 	submanifold of $(\overline{M},\overline{g})$ with $r>1.$
 	Suppose $\mathcal{U}$ is a coordinate neighbourhood of $M$ 
 	and $\{\xi_i\}; i\in\{1,\cdots, r\}$  a basis of
 	$\Gamma(RadTM)\arrowvert_{\mathcal{U}}$. Then there exist smooth sections
 	$\{N_i\}$ of $\mathcal{S}(TM^\bot)^\bot\arrowvert_{\mathcal{U}}$ such that
 	
 	\begin{equation}\label{nultrans gene}
 	\overline{g}(N_i, \xi_i)=\delta_{ij}, 
 	\end{equation}
 	and 
 	\begin{equation}\label{transnull}
 	\overline{g}(N_i, N_j)=0
 	\end{equation}
 	for any $i,j\in\{1, ..., r\}.$
 \end{theorem}
 
 It follows that there exists a null transversal vector bundle $ltr(T M)$ 
 locally spanned by $\{Ni\}$. By using (\ref{nultrans gene}) and (\ref{nullscreen}) it is
 easy to check that $B=\{\xi_1,\cdots,\xi_r, N_1,\cdots,N_r\}$ is a basis of
 $\Gamma{S}(TM^\bot)^\bot_{|\mathcal{U}}$. 
 The set of local sections $\{N_i\}$ is not unique enven if one use the same vector bundle  in general.
 Let $tr(TM)$ be complementary called (but not orthogonal) vector 
 bundle to $T M$ in $T \overline{M}\arrowvert_M$. Then the following hold
 \begin{equation}\label{I6}
 tr(T M) = ltr(T M)\oplus_{Orth} \mathcal{S}(T M^\perp ),
 \end{equation}
 \begin{equation}\label{tranvdecomp}
 T\overline{ M}\arrowvert_M  
 = \mathcal{S}(T M)\oplus_{Orth} (RadT M \oplus ltr(T M))\oplus_{Orth}\mathcal{S}(T M^\perp)
 = T M \oplus tr(T M).
 \end{equation}
 
 Let $\overline{\nabla}$ be the Levi-Civita connection on $\overline{M}$. 
 As $T M$  and $tr(TM )$ are complementary sub-bundles of $T\overline  M\arrowvert_M$,
 the Gauss and Weingarten formulae are
 \begin{equation}\label{IGauss}
 \overline{\nabla}_XY=\nabla_XY+h(X,Y),
 \end{equation}
 \begin{equation}\label{Iweigart}
 \overline{\nabla}_XV=-A_VX+\nabla_X^tV,
 \end{equation}
 $\forall X,Y\in \Gamma(TM), V\in \Gamma(tr(TM))$. $\nabla$ and $\nabla^t$ are linear
 connections on $TM$ and the vector bundle $tr(TM)$ called  the induced linear
 connection and the transversal linear connection on $M$ respectively.
 
 $\nabla$  is torsion-free linear connection. The components $\nabla_XY$ 
 and $-A_VX$ belong to $\Gamma(TM)$,
 $h(X,Y)$ and $\nabla_X^tV$  to $\Gamma(tr(TM))$.
 Also $h$ is a $\Gamma(tr(TM))$-valued symmetric bilinear form on $\Gamma(TM)$  
 called the second fundamental form of $M$ with respect to $tr(TM)$.  
 $A$ is a $\Gamma(TM)$-valued bilinear form defined on $\Gamma(tr(TM))\times \Gamma(TM)$
 called shape operator of $M$.
 From the geometry of Riemannian submanifolds \cite{Chen} and non-degenerate 
 submanifolds \cite{ONei}, it is known that the induced connection on a non-degenerate
 submanifold is a Levi-Civita connection.
 Unfortunately, in general, this is not true for a null submanifold.
 
 \begin{equation}\label{I23}
 (\nabla_Xg)(X, Y)=\overline{g}(h^l(X, Y)+\overline{g}(h^l(X,Z), Y)
 \end{equation}
 
 \begin{equation}\label{I24}
 (\nabla_X^t\overline{g})(V, V^{'})=-\{\overline{g}(A_VX, V^{'})+\overline{g}(A_{V^{'}}X, V)\}
 \end{equation}
 
 for all $X, Y, Z \in \Gamma(TM);, V, V^{'}\in \Gamma(tr(TM))$. Thus, it follows 
 that the induced connection $\nabla$ is not a Levi-Civita connection.
 
 According to the decomposition (\ref{I6}), let $L$ and $S$ denote the projection morphisms of $tr(TM)$
 onto $ltr(TM)$ and $\mathcal{S}(TM^\perp)$ respectively, $h^{l} = L\circ h$, $h^{s} = S\circ h$, 
 $D^{l}_{X}V= L(\nabla_X^tV)$, $D_X^sV= S(\nabla_X^tV)$. The transformations 
 $D^{l}$ and $D^{s}$ do not define linear connections but Otsuki connections on $tr(TM)$ with respect to the
 vector bundle morphisms $L$ and $S$.  Then we have  
 \begin{equation}\label{I10}
 \overline{\nabla}_XY=\nabla_XY+h^l(X, Y)+h^s(X, Y)
 \end{equation}
 (\ref{nultrans gene}).
 \begin{equation}\label{I11}
 \overline{\nabla}_XN=-A_NX+D_X^lN+D^s(X, N)
 \end{equation}
 
 \begin{equation}\label{I12}
 \overline{\nabla}_XW=-A_WX+D^l(X, W)+\nabla_X^sW
 \end{equation}
 $\forall X,Y\in \Gamma(TM),N\in \Gamma(ltr(M))$ and $W\in \Gamma(S(TM^\perp)).$
 Since $\overline{\nabla}$ is a metric connection, using (\ref{I10})-(\ref{I12}) we have 
 \begin{equation}\label{I13}
 \overline{g}(h^s(X, Y),W)+\overline{g}(Y,D^l(X, W))=g(A_WX, Y)
 \end{equation}
 \begin{equation}\label{I14}
 \overline{g}(D^s(X, N), W)=\overline{g}(N, A_WX).
 \end{equation}
 As $h^l$ and $h^s$ are $\Gamma(ltr(TM))$-valued and 
 $\Gamma(\mathcal{S}(TM^\bot))$-valued respectively, we call
 them the null second fundamental form and the screen second fundamental form of $M$.
 
 Suppose $M$ is either r-null with $r<min\{m, k\}$
 or coisotropic. Then, using the decomposition (\ref{Radical}) we get
 
 \begin{equation}\label{gauss ecran}
 \nabla_XY=\nabla_X^{*}PY+h^{*}(X, PY)  
 \end{equation}
 
 \begin{equation}\label{I18}
 \nabla_X\xi=-A_\xi^{*} X+\nabla_X^{*t}\xi
 \end{equation}
 
 for any $X, Y\in \Gamma(TM)$ and $\xi\in \Gamma(RadTM)$, where $\{ \nabla_X^{*}PY, 
 A_\xi^{*} X\}$
 and $\{h^{*}(X, PY),\nabla_X^{*t}\xi\}$ belong to $\Gamma(\mathcal{S}(TM))$ 
 and $\Gamma(Rad(TM))$ respectively. It follows
 that $\nabla^{*}$ and $\nabla^{*t}$ are linear connections on $\mathcal{S}(TM)$ 
 and $RadTM$ respectively. On the other hand, $h^{*}$ 
 and $A^{*}$ are $\Gamma(Rad T M )$-valued and $\Gamma(\mathcal{S}(T M ))$-valued bilinear 
 forms on $\Gamma(T M ) \times \Gamma(\mathcal{S}(T M ))$
 and $\Gamma(Rad T M ) \times \Gamma(T M )$ Called the second fundamental forms of 
 distributions $\mathcal{S}(T M )$ and $Rad(T M )$ respectively.\\
 For any $\xi\in\Gamma(RadTM)$ consider the linear operator
 $$A_\xi:\Gamma(TM)\rightarrow\Gamma(\mathcal{S}(TM));A_\xi^{*} X=A^{*}(\xi,X),
 ,\forall,X\in \Gamma(TM)$$
 and call it the shape operator of $\mathcal{S}(TM)$ with respect to $\xi$. Also, 
 call  $\nabla^{*}$ and $\nabla^{*t}$ the induced connections on $\mathcal{S}(TM)$
 and $RadTM$ respectively.
 
 The second fundamental form and shape operator of non-degenerate submanifold of a semi-Riemannian manifold
 are related by means of the metric tensor field. Contrary in case of null submanifolds
 there are interrelations between geometric objects induced by $tr(TM)$ on one side 
 and geometric objects induced by $S(TM)$ on the other side. More precisely,

 \begin{equation}\label{form}
 \overline{g}(h^l(X, PY)=g(A_\xi^{*} X, PY),\ \overline{g}(h^{*}(X, PY),N)=g(A_NX, PY),
 \
 \overline{g}( h^l(X, \xi), \xi)=0, A_\xi^{*} \xi=0.
 \end{equation}

 From (\ref{form}) as $h^l$ is symetric, it follows that the shape operator of 
 $\mathcal{S}(TM)$ is a self-adjoint operator on $\mathcal{S}(TM)$.
 
 Next, consider a coordinate neighbourhood $\mathcal{U}$ of $M$ and let $\{N_i , W_{\alpha} \}$ be a
 basis of $\Gamma(tr(T M )\arrowvert_M )$ where $N_i \in \Gamma(ltr(T M )\arrowvert_M ), i \in\{1, . . . , r\}$
 and $W_{\alpha} \in \Gamma(\mathcal{S}(T M^\bot)\arrowvert_{\mathcal{U} }), \alpha \in \{r + 1, . . . , k\}$. 
 Then ( \ref{I10}) becomes
 
 \begin{equation}\label{I22}
 \overline{\nabla}_XY=\nabla_XY+\sum_{i=1}^r h_i^l(X, Y)N_i+\sum_{\alpha=r}^k h_i^s(X, Y)W_{\alpha}.
 \end{equation}
 We call $\{h_i^l\}$ and $\{h^s_{\alpha}\}$ the local null second fundamental forms and 
 the local screen second fundamental forms of $M$ on $\mathcal{U}$.
 
 We recall  the equations of Gauss, Codazzi and Ricci which play an important role in studying 
 differential geometry of non-degenerate submanifolds(See \cite{Duggal} for details).
 Let $\overline{R}$ and $R$ denote the Riemannian curvature tensors 
 of $\overline{\nabla}$ and $\nabla$ on $\overline{M}$ and $M$ respectively. The Gauss
 equation are given by
 \begin{eqnarray}\label{I25}
 \overline{R}(X, Y)Z &=& R(X, Y)Z+A_{h^l(X, Z)}Y-A_{h^l(Y, Z)}X+A_{h^s(X, Z)}Y-A_{h^s(Y, Z)}X\nonumber\\
 & &+(\nabla_Xh^l)(Y, Z)-(\nabla_Yh^l)(X, Z)+D^l(X, h^s(Y, Z))\nonumber\\
 &&-D^l(Y,h^s(X, Z))
 +(\nabla_Xh^s)(Y, Z)-(\nabla_Yh^s)(X, Z)\nonumber\\
 &&+D^s(X,h^l(Y, Z))-D^s(Y,h^s(X, Z))\nonumber
 \end{eqnarray}
 $\forall X,Y ,Z, U\in \Gamma(TM).$
 Therefore
 \begin{eqnarray}\label{I26}
 \overline{R}(X, Y, PZ, PU)&=&R(X, Y, Z, PU)+\overline{g}(h^{*}(Y, PU),h^l(X, Z))\nonumber\\
 &&-\overline{g}(h^{*}(X, PU), h^l(Y, Z))
 +\overline{g}(h^s(Y, PU),h^s(X, Z))\nonumber\\
 &&-\overline{g}(h^s(X, PU),h^s(Y, Z)),
 \end{eqnarray}
 
 Throughout, we consider that the submanifold is coisotropic that is $h^s=0$.
 Let the differential $1-$forms $\omega_i,i=1,\cdots,r (0<r<min\{m,n\})$ defined by:
 \begin{equation}\label{oneforma}   
 \omega_i(X)=\overline{g}(X, N_i),, \forall, X\in \Gamma(TM\arrowvert_{\mathcal{U}}),i\in \{1, ..., r\}.
 \end{equation} 
 Then any vector $X$ on $M$ is expressed on $\mathcal{U}$ as follows
 \begin{equation}\label{projectionecran}
 X=PX+\sum_{i=1}^r\omega_i(X)\xi_i
 \end{equation}
 where $P$ is the projection morphism of $\Gamma((TM))$ onto $\Gamma(S(TM))$.
 
 
 \section{Normalization and Induced semi-Riemannian metric}\label{Section3}
 Let $\varpi$ denote the $1$-form wich satisfying $\varpi(.) = \overline{g}(N,.)$.
 Then, take
 \begin{equation}
 \omega = f^\star{\varpi}
 \end{equation}
 to be its restriction to $M$ and the map $f : M\rightarrow\overline{M} $ the inclusion map.
 Throughout, a screen distribution on $M$ is denote by $\mathcal{S}(N)$.
 We define the associated semi-Riemannian metric on $\overline{M}$ as
 \begin{equation}
 \underline{g}=\overline{g}-\varpi\otimes\varpi\quad\mbox{and}\quad
 \widetilde{ g}=i^{\star}\underline{g}
 =g-\omega_{i=1}^r\otimes\omega_{i=1}^r. \label{garc}%
 \end{equation} 
 \begin{lemma}\label{lem1}
 	Let $(M^n,g)$ be a r-null submanifold in semi-Riemannian manifold 
 	($\overline{M}^{n+1},\overline{g})$.
 	Then, $\widetilde{ g}$ is  nondegenerate.
 \end{lemma}
 \begin{prof7}
 	Take $V_x=\sum_{i=1}^r\xi_{{i}_|x}$, 
 	then  
 	$ \sum_{i=1}^r\widetilde{ g}(U_x,\xi_{{i}_|x})=
 	0= \sum_{i=1}^r g(u,\xi_{{i}_|x})-\sum_{i=1}^r\omega_i(U)\omega_i(\xi_{{i}_|x})$,
 	$\ \sum_{i=1}^r \omega_i(U)=0$. 
 	Then
 	\begin{equation}\label{nondegenerate1}
 	U\in \mathcal{S}(N).
 	\end{equation} 
 	\begin{equation}\label{nondegenerate2}
 	0=\widetilde{ g}(U,V)=g(U,V) \Rightarrow U\in T_xM^\bot.
 	\end{equation}
 	From $(1)$ and $(2)\Rightarrow U\in \mathcal{S}(N_x)\cap T_xM^\bot =\{0\}$ and then 
 	$U=0$ which show $\widetilde{ g}$ is nondegenerate.
 \end{prof7}
 
 With $\sum_{i=1}^r\omega_i=\overline{g}(N_i,.)$.
 It is nothing to see that 
 \begin{equation}
 \widetilde{ g}(\xi_i,\xi_j)=-\delta_{ij}, \omega_i(.)=\widetilde{ g}(\xi_i,.).
 \end{equation} 
 \begin{lemma}
 	Let $(M,g,\mathcal{S}(N))$ be a r-null submanifold in semi-Riemannian manifold, then
 	$(M,\widetilde{ g})$ is a semi-Riemannian manifold of index$~ q-1$. In particular,
 	$(M,\widetilde{ g})$ is Riemannian if
 	$(\overline{M},\overline{g})$ is Lorentzian and $(M,\widetilde{g})$ is Lorentzian if 
 	$(\overline{M},\overline{g})$ is of index 2.
 \end{lemma}
 \begin{prof7}
 	Take  $(e_{2},\dots,e_{n+1})$ an orthonormal basis of
 	$\mathcal{S}(N)$ for $\overline{g}$ and note that 
 	$(\xi,e_{2},\dots,e_{n+1})$ is an orthonormal basis of $TM$ for the rigged metic. Elsewhere, it is a fact that
 	$\mathcal{S}(N)$ is a non degenerate distribution of index $q-1$ (see \cite{Duggal}). 
 	Since $\widetilde{ g}$ agrees with ${g}$ on $\mathcal{S}(N)$ and $\widetilde{ g}(\xi,\xi)=-\delta_{ij}$, it follows that 
 	$(M,\widetilde{ g})$ is a semi-Riemannian manifold of index $~q-1$.
 \end{prof7}
 \begin{lemma}\label{lemme} 
 	For all $ X,Y,Z\in\Gamma(TM)$, we have 
 	\begin{eqnarray}\label{equa11}
 	(\nabla _X \widetilde g)(Y,Z)
 	&=&\sum_{i=1}^r\omega_i(Y)[g(A_{\xi_i}^{\star} X, PZ)- g(A_{N_i}X,PZ)]\nonumber\\
 	&&+ \sum_{i=1}^r\omega_i(Z)[ g(A_{\xi_i}^{\star} X,PY)- g(A_{N_i}X,PY)]\nonumber\\
 	&&+2\sum_{i=1}^r\tau^{N_i}(X)\omega_i(Y)\omega_i(Z).\nonumber\\
 	\end{eqnarray}
 \end{lemma}
 Let $(M,g,\mathcal{S}(N))$ be a normalized null submanifold of a semi-Riemannian manifold,
 then the integral curves of the rigged vector field $\xi$ are pregeodesic but not
 geodesic in general. The following lemma shows that in case the normalization
 is a conformal vector field, then $\xi$ is $\overline{g}$-geodesic.
 \begin{lemma}
 	\label{rig con}
 	Let $ (M,g,\mathcal{S}(N)) $ be a conformal normalized  null submanifold  of a semi-Riemannian
 	manifold such, then $\xi_i$ is $\overline{g}$-geodesic and $ \sum_{i=1}^r\tau(\xi_i) = 0$.
 \end{lemma}
 \begin{prof7}
 	Firstly, note that  $ \sum_{i=1}^r\tau(\xi_i)= 
 	\sum_{i=1}^r\overline{g}(\overline{\nabla}_{\xi_i} N_i, \xi_i)$
 	$N_i$ being conformal, there exists a function  $\lambda$ such that  
 	$\overline{g}(\overline{\nabla}_XN_i, Y) + \overline{g}(\overline{\nabla}_YN_i, X) 
 	= \lambda\overline{g}(X,Y).$
 	Take $ X = Y =  \sum_{i=1}^r\xi_i$ to get  
 	$ 2 \sum_{i=1}^r\overline{g}(\overline{\nabla}_{\xi_i}N_i, \xi_i)
 	=  \sum_{i=1}^r\lambda\overline{g}(\xi_i, \xi_i) = 0;$ hence 
 	$ \sum_{i=1}^r\tau(\xi_i)
 	=  \sum_{i=1}^r\overline{g}(\overline{\nabla}_{\xi_i}N_i, \xi_i) 
 	= 0$. 
 	It follows that $ \sum_{i=1}^r\overline{\nabla}_{\xi_i}\xi_i
 	= - \sum_{i=1}^r\tau(\xi_i)\xi_i = 0.$
 \end{prof7}
 
 \begin{definition}\label{conformals}
 	A normalized null submanifold $(M,g,N)$ of a semi-Riemannian manifold
 	$ (\overline{M}, \overline{g}) $ is said to have a conformal  screen~\cite{ATD} 
 	if there exists a non vanishing smooth function $\varphi$ on $M$ such that 
 	$ A_N = \varphi A_{\xi}^{\star}$ holds.
 \end{definition}
 This is equivalent to saying that  $ g(A_NX,PY) = \varphi g(A_\xi^{\star} X,Y)$ for all tangent 
 vector fields $ X $ and $ Y $. The function $\varphi $ is called the conformal factor.
 \section{ Relation between the null and the associated semi-Riemannian Geometry} \label{Section4}
 
 The main focus of this section lies on deriving jump formulas for the various curvature
 quantities, that is, how the Riemann and Ricci tensor and scalar curvature
 of course the reason why this is of a particular interest lies in physics, mainly general 
 relativity,where such formulas might find applications due to the Einstein field equations.
 In this respect, using the lemma (\ref{lemme})
 and as  connections $\widetilde{\nabla}$ and $\nabla$ are torsion 
 free, we can first prove the relation between  the Levi-Civita connection of 
 the non-degenerate associate metric $ \widetilde{ g}$ and the rigged connection $\nabla$ 
 in the following.
 \begin{proposition} \label{connectioneta}
 	Let $(M^n,g)$ be a  $r-$normalized coisotropic null submanifold in semi-Riemannian Manifold with rigged 
 	vector field $\xi$.
 	Then, for all $ X,Y \in \Gamma(TM)$, we have
 	\begin{eqnarray}\label{equa15}
 	\widetilde{\nabla}_{ X}Y&=&\nabla_{ X}Y    
 	+\frac{1}{2}\sum_{i=1}^r\big[ 2 \langle A_{\xi_i}^{*} X,Y\rangle
 	-\langle A_{N_i}X,Y\rangle-\langle A_{N_i}Y,X\rangle\nonumber\\
 	&&+\sum_{i=1}^r\omega_i(X)\tau^{N_i}(Y)
 	+\sum_{i=1}^r\omega_i(Y)\tau^{N_i}(X)\big]\xi_i        
 	+\sum_{i=1}^r\omega_i(X)(i_{Y}d\omega_i)^{\sharp_{\omega_i}}
 	+\sum_{i=1}^r\omega_i(Y)(i_{X}d\omega_i)^{\sharp_{\omega_i}}.                                    \nonumber 
 	\end{eqnarray}
 	
 	In particular for a closed normalization,
 	\begin{eqnarray}\label{equa16}
 	\widetilde{\nabla}_{ X}Y&=&\nabla_{X}Y    
 	+\frac{1}{2}\sum_{i=1}^r\big[ 2 \langle A_{\xi_i}^{*} X,Y\rangle-
 	\langle A_{N_i}X,Y\rangle-\langle A_{N_i}Y,X\rangle\big. \nonumber\\
 	&&+\big.\sum_{i=1}^r\omega_i(X)\tau^{N_i}(Y)
 	+\sum_{i=1}^r\omega_i(Y)\tau^{N_i}(X)\big]\xi_i.\nonumber\\              \end{eqnarray}
 \end{proposition}
 Let $\widetilde{R}$ and $R$ denote the Riemann curvature tensors of
 $\widetilde{\nabla}_{\alpha}$
 and $\nabla$ respectively.
 Using proposition (\ref{connectioneta}) we prove the following. 
 \begin{proposition}\label{tensoreta}
 	Let $(M, g^n, N )$ be a $r-$closed normalized coisotropic null submanifold of semi-Riemannian Manifold with
 	rigged vector field $\xi$ . Then, for all $X, Y, W \in \Gamma(T M )$ and $U \in T M^\bot$ we have,
 	\begin{eqnarray}\label{prof1}
 	\widetilde{ g}(\widetilde{R}(X,Y)Z,PW) 
 	&=& \langle R(X,Y)Z,PW\rangle\nonumber\\
 	&&+\frac{1}{2} 
 	\sum_{i=1}^r\{ 2 g(A_{\xi}^{\star}X,Z) 
 	-\sum_{i=1}^r\langle A_{N_i}X,Z\rangle
 	- \sum_{i=1}^r\langle A_{N_i}Z,X\rangle \nonumber\\
 	&&+\tau^{N_i}(X)\sum_{i=1}^r\omega_i(Z)+\tau^{N_i}(Z)
 	\sum_{i=1}^r\omega_i(X)\}g(A_\xi^{\star}(Y,W) \nonumber\\ 
 	&&-\{2 g(A_\xi^{\star} Y,Z)- \sum_{i=1}^r\langle A_{N_i}Y,Z\rangle
 	- \sum_{i=1}^r\langle A_{N_i}Z,Y\rangle \nonumber\\
 	&&+\tau^{N_i}(Y)\sum_{i=1}^r\omega_i(Z)+\tau^{N_i}(Z)
 	\sum_{i=1}^r\omega_i(Y)\} g(A_\xi^{\star} X,W)\},
 	\end{eqnarray}
 	\begin{eqnarray}\label{prof2}
 	\widetilde{ g}(\widetilde{R}(X,Y)Z,U) 
 	&=&- \langle R(X,Y)U,PZ\rangle \nonumber\\
 	&&-\frac{1}{2} \sum_{i=1}^r\big[ \langle A_{N_i}U,Y\rangle 
 	g(A_\xi^{\star} X,Z)-\langle A_{N_i}U,X\rangle g(A_\xi^{\star} Y,Z)\big ]\nonumber\\
 	&&-\frac{1}{2} \sum_{i=1}^r\big[\tau^{N_i}(X)
 	g(A_\xi^{\star} Y,Z)-\tau_i^N (Y)g(A_\xi^{\star} X,Z)\big]\omega_i(U),
 	\end{eqnarray}
 \end{proposition}

 \end{document}